\RequirePackage{fix-cm}
\documentclass[condensed]{svjour3}       
\smartqed  
\usepackage{graphicx}
\usepackage{amssymb}
\usepackage{amsmath}
\newcommand{\ds}{\displaystyle}

\begin{document}

\title{Analysis of viscoelastic flow with a generalized memory and its exponential convergence to steady state\thanks{The work was supported by the Natural Science Foundation of Shaanxi Province under grant No. 2019JQ-398, by the Scientific Research Program of Shaanxi Provincial Education Department under grant No. 19JK0841, by the China Postdoctoral Science Foundation under grants No. 246135 and No. 2021TQ0017, and by the International Postdoctoral Exchange Fellowship Program (Talent-Introduction Program) No. YJ20210019.}
}

\titlerunning{Exponential convergence of viscoelastic flow with memory}        

\author{Yingwen Guo        \and
        Xiangcheng Zheng$^*$ 
}


\institute{Y. Guo \at
             School of Mathematical Sciences, Peking University, Beijing 100871, China. 
           \and
           X. Zheng \at
              School of Mathematical Sciences, Peking University, Beijing 100871, China.\\
              \email{zhengxch@math.pku.edu.cn}          
}

\date{Received: date / Accepted: date}

\maketitle

\begin{abstract}
We investigate a viscoelastic flow model with a generalized memory, in which a weak-singular component is introduced in the exponential convolution kernel of classical viscoelastic flow equations that remains untreated in the literature. We prove the well-posedness and regularity of the solutions, based on which we prove the exponential convergence of the solutions to the steady state. The proposed model serves as an extension of classical viscoelastic flow equations by adding a dimension characterized by the power of the weak-singular kernel, and the derived results provide theoretical supports for designing numerical methods for both the considered equation and its steady state.
\end{abstract}

\section{Introduction}
Viscoelastic flows arise in widely applications and have attracted extensive attentions. A typical governing equation is derived from the combination of the Oldroyd's model, the continuity equation and the momentum conservation equation as follows
\begin{equation}\label{mh0}
\mathbf{u}_t-\mu\Delta \mathbf{u}+(\mathbf{u}\cdot\nabla)\mathbf{u}-
\rho\int_0^te^{-\delta(t-s)}\Delta \mathbf{u}ds+\nabla
p=\mathbf{f}.
\end{equation}
Here $\mathbf{u}=\mathbf{u}(\mathbf{x},t)=(u_1(\mathbf{x},t),
u_2(\mathbf{x},t))$ is the velocity
vector, $p=p(\mathbf{x},t)$ is the pressure of the fluid, $\mathbf{f}=(f_1(\mathbf{x},t),f_2(\mathbf{x},t))$ is the prescribed external
force, $\mathbf{x}\in\Omega$ for some two-dimensional convex polygonal domain $\Omega$ with the boundary $\partial \Omega$, $\mu$ is the solvent viscosity, $\delta=1/\lambda_1$ where $\lambda_1$ refers to the relaxation time and $\rho=(\mu/\lambda_1)(\lambda_1/\lambda_2-1)>0$ where $\lambda_2$ stands for the retardation time satisfying the restriction
$0\leq\lambda_2\leq\lambda_1$.

Extensive mathematical and numerical analysis for model (\ref{mh0}) and related problems have been carried out, see e.g., \cite{AraMen,CanEwi,Coc1,Coc2,Gal,LiHe,LiShe,Lia,Lub,MusMcl}. In particular, the convergence of (\ref{mh0}) to a steady state in an exponential rate with respect to time was rigorously proved in \cite{HeLin,Kot,Ort,WanLin}. A main application of this result lies in approximating the steady state of (\ref{mh0}) with small viscosity $\mu$ via computing (\ref{mh0}) \cite{He2013}. To be specific, solving the steady state of (\ref{mh0}) requires certain iterative method due to its strong nonlinear term $(\mathbf{u}\cdot\nabla)\mathbf{u}$. However, it was shown in \cite{He2013,HeLi} that the uniqueness condition of the Stokes iterative, which is an efficient iterative method with a constant
coefficient matrix compared with other iterative methods, may be lost for small viscosity $\mu$. Thus one may alternatively solve (\ref{mh0}) to approximate its steady state, which requires the convergence estimates of (\ref{mh0}) to its steady state.

 Compared with classical Naiver-Stokes equations, which corresponds to $\rho=0$ in (\ref{mh0}), a mainly overcomed difficulty in the investigation of model (\ref{mh0}) in the literature stems from the newly encountered convolution term, which introduces nonlocal features in time by an exponential kernel to model the memory effects of the fluid and thus the dynamics. In the past decades, increasingly experimental evidences show that the power function kernels in, e.g., the fractional calculus, adequately describe the creep and relaxation of viscoelastic materials \cite{BagTor,Koe,Mer,Mai} and thus provide more accurate modeling for the memory effects in the dynamics of viscoelastic fluids \cite{Ali2,Kha,Pao,PerKar,Suz}. Nevertheless, the corresponding viscoelastic flow equation like (\ref{mh0}) involving the power function memory kernel is meagerly studied in the literature.

Motivated by these discussions, we consider the following viscoelastic flow model with a generalized memory kernel combining the exponential component as model (\ref{mh0}) and the newly added power-law decay factor $t^{-\beta}$ for $0\leq\beta<1$
\begin{eqnarray}
\qquad \left\{
\begin{array}{@{}ll}
\ds\mathbf{u}_t-\mu\Delta \mathbf{u}+(\mathbf{u}\cdot\nabla)\mathbf{u}-
\rho\int_0^t(t-s)^{-\beta}e^{-\delta(t-s)}\Delta \mathbf{u}ds+\nabla
p=\mathbf{f},\\[0.05in]
(\mathbf{x},t)\in\Omega\times [0,\infty);~~\mbox{div}\;\mathbf{u}=0,~(\mathbf{x},t)\in \Omega\times [0,\infty);\\[0.05in]
\mathbf{u}(\mathbf{x},0)=\mathbf{u}_0(\mathbf{x}),~\mathbf{x}\in \Omega;~~ ~\mathbf{u}(\mathbf{x},t)|_{\partial\Omega}=0,~~ t\in
[0,\infty).
\end{array}
\right. \label{eq3.1}
\end{eqnarray}
The combined memory kernel $t^{-\beta}e^{-\delta t}$ appears in many applications, see e.g., \cite[Equation 1.4]{Xu} and \cite[Equation 1.22]{Lin}. Compared with (\ref{mh0}), the additional factor $t^{-\beta}$ introduces initial singularities and the convolution kernel is no longer continuously differentiable as that in (\ref{mh0}). For this reason, the analysis needs careful treatments near the initial time. Furthermore, the indefinite integral of this combined kernel could not be evaluated into a closed form as its components $e^{\delta t}$ and $t^{-\beta}$, which significantly complicates the analysis.

In this paper, we address the aforementioned concerns to analyze model (\ref{eq3.1}) and prove its exponential convergence to the following steady state
\begin{eqnarray}
\qquad \left\{
\begin{array}{@{}ll}
\ds-\bigg(\mu+\frac{\rho\Gamma(1-\beta)}{\delta^{1-\beta}}\bigg)\Delta \bar{\mathbf{u}}+(\bar{\mathbf{u}}\cdot\nabla)\bar{\mathbf{u}}+\nabla
\bar{p}=\bar{\mathbf{f}},~~\mathbf{x}\in \Omega,\\[0.15in]
\mbox{div}\;\bar{\mathbf{u}}=0,~\mathbf{x}\in \Omega;~~\bar{\mathbf{u}}(\mathbf{x})|_{\partial\Omega}=0
\end{array}
\right. \label{eq1.1}
\end{eqnarray}
which provide theoretical supports for designing numerical methods for both the considered equation (\ref{eq3.1}) and its steady state (\ref{eq1.1}). Note that when $\beta$ tends to $0$, this equation approaches to the steady state of the traditional viscoelastic flow equation (\ref{mh0}) proposed in, e.g., \cite[Equation 1.3]{HeLin}. That is, the current work mathematically generalizes the existing models and results by introducing an extra dimension characterized by the parameter $\beta$.

The rest of the paper is organized as follows: In Section \ref{sec2} we introduce notations and preliminary results to be used subsequently. In Section \ref{sec3} we prove well-posedness and regularity results for model (\ref{eq3.1}), based on which we prove exponential convergence of viscoelastic flow equation (\ref{eq3.1}) to its steady state (\ref{eq1.1}) in Sections \ref{sec4}-\ref{sec5}.
\section{Preliminaries}\label{sec2}
\subsection{Spaces and notations}
We follow \cite{adams,girault,heywood1,te1} to introduce useful spaces and notations.
Define the following Hilbert spaces
$$
\mathbf{X}=H^1_0(\Omega)^2,~~ \mathbf{Y}=L^2(\Omega)^2,~~
M=L^2_0(\Omega)=\left\{q\in L^2(\Omega);  \int_{\Omega}qd\mathbf{x}=0\right\},$$
where $L^2(\Omega)^d$~($d=1,2$) is equipped
with the usual $L^2$-scalar product $(\cdot,\cdot)$ and $L^2$-norm
 $\|\cdot\|_{0}$, and $H^1_0(\Omega)$ and $\mathbf{X}$ are equipped with the following scalar product and equivalent norm
$$
(\nabla \cdot,\nabla \cdot),~~~|\cdot|_1=\|\nabla \cdot\|_{0}.
$$
Here, $\|\cdot\|_i$ and $|\cdot|_i$ denote the usual norm and semi norm of the sobolev space $H^i(\Omega)^d$, for $i=0,1,2$. Define the closed subset $\mathbf{V}$ of $\mathbf{X}$ and the closed subset $\mathbf{H}$ of $\mathbf{Y}$ as follows
$$
\mathbf{V}=\{\mathbf{v}\in \mathbf{X}; {\rm div} \;\mathbf{v}=0\},~~\mathbf{H}=\{\mathbf{v}\in \mathbf{Y}; {\rm div} \;\mathbf{v}=0, \mathbf{v}\cdot \mathbf{n}|_{\partial\Omega}=0\}.
$$
Furthermore, let $P$ be the $L^2$-orthogonal projection of $\mathbf{Y}$ onto $\mathbf{H}$ and $A$ be the Stokes operator defined by $A=-P\Delta$ with $\mathbf{D}(A):=H^2(\Omega)^2\cap
\mathbf{V}$, which satisfies
 \cite{adams,He2013,heywood1,larsson}
\begin{align}
&\|\mathbf{v}\|_{2}^2\leq c\|A\mathbf{v}\|_{0}^2,~~ \mathbf{v}\in \mathbf{D}(A),
\label{eq2.3}\\
\gamma_0\|\mathbf{v}\|_{0}^2&\leq |\mathbf{v}|_{1}^2,~~ \mathbf{v}\in
\mathbf{X},~~\gamma_0|\mathbf{v}|_{1}^2\leq \|A\mathbf{v}\|_{0}^2,~~ \mathbf{v}\in \mathbf{D}(A),\label{eq2.4}
\end{align}
where $\gamma_0$, $c>0$ are generic constants
depending on $\Omega$. In particular, $A^{\frac{1}{2}}$ satisfies $\|A^{\frac{1}{2}}\mathbf{v}\|_0=|\mathbf{v}|_1$ for $\mathbf{v}\in \mathbf{V}$.

With the above definitions, the bilinear forms
$a(\cdot,\cdot)$ and $d(\cdot,\cdot)$ on $\mathbf{X}\times \mathbf{X}$ and $\mathbf{X}\times
M$ are respectively defined by
$$
a(\mathbf{u},\mathbf{v})=(\nabla \mathbf{u},\nabla \mathbf{v}),~~d(\mathbf{v},q)=-(\mathbf{v},\nabla q)=(q,{\rm div} \;\mathbf{v}),~~ \mathbf{u},\mathbf{v}\in \mathbf{X}, ~~q\in M.
$$
The bilinear form $a(\cdot,\cdot)$ is continuous and coercive on $\mathbf{X}\times \mathbf{X}$, and $d(\cdot,\cdot)$ is continuous on $\mathbf{X}\times M$ and satisfies the well-known inf-sup condition: there exist a positive constant $C_0>0$ such that for all $q\in M$
\begin{eqnarray}
\sup_{\mathbf{v}\in \mathbf{X}}\frac{|d(\mathbf{v},q)|}{|\mathbf{v}|_1}\geq C_0\|q\|_0.\label{eq5.1}
\end{eqnarray}

In addition, the trilinear form $a_1(\cdot,\cdot,\cdot)$ on $\mathbf{X}\times \mathbf{X}\times \mathbf{X}$ is defined by
\begin{eqnarray}
a_1(\mathbf{u},\mathbf{v},\mathbf{w})=\langle(\mathbf{u}\cdot\nabla)\mathbf{v},\mathbf{w}\rangle_{\mathbf{X}',\mathbf{X}},\ \ \ \mathbf{u},\mathbf{v},\mathbf{w}\in \mathbf{X},\nonumber
\end{eqnarray}
which is continuous on $\mathbf{X}\times \mathbf{X}\times \mathbf{X}$ and satisfies \cite{heywood1,te1}
\begin{align}
a_1(\mathbf{u},\mathbf{v},\mathbf{w})=&-a_1(\mathbf{u},\mathbf{w},\mathbf{v}),\ \ \ \forall \ \mathbf{u}\in \mathbf{V},\mathbf{v},\mathbf{w}\in \mathbf{X},\label{eq2.01}\\
|a_1(\mathbf{u},\mathbf{v},\mathbf{w})|\leq & N|\mathbf{u}|_1|\mathbf{v}|_1|\mathbf{w}|_1,\ \ \ \forall\  \mathbf{u},\mathbf{v},\mathbf{w}\in \mathbf{X},\label{eq2.02}\\
|a_1(\mathbf{u},\mathbf{v},\mathbf{w})|\leq & c_0\|\mathbf{u}\|_0^{\frac{1}{2}}|\mathbf{u}|_1^{\frac{1}{2}}|\mathbf{v}|_1^{\frac{1}{2}}\|A\mathbf{v}\|_0^{\frac{1}{2}}\|\mathbf{w}\|_0,\ \ \ \forall\ \mathbf{u}\in \mathbf{X}, \mathbf{v}\in \mathbf{D}(A), \mathbf{w}\in \mathbf{Y},\label{eq2.03}\\
|a_1(\mathbf{v},\mathbf{u},\mathbf{w})|\leq & c_0\|\mathbf{v}\|_0^{\frac{1}{2}}\|A\mathbf{v}\|_0^{\frac{1}{2}}|\mathbf{u}|_1\|\mathbf{w}\|_0,\ \ \ \forall\ \mathbf{\mathbf{u}\in \mathbf{X}}, \mathbf{v}\in \mathbf{D}(A), \mathbf{\mathbf{w}\in \mathbf{Y},}\label{eq2.04}
\end{align}
where $c_0$ is a positive constant depending only on $\Omega$ and $N$ is defined in terms of
\begin{align}
b(\mathbf{u},\mathbf{v},\mathbf{w}):&=((\mathbf{u}\cdot\nabla)\mathbf{v},\mathbf{w})+\frac{1}{2}(({\rm div} \mathbf{u})\mathbf{v},\mathbf{w})\nonumber\\
&\ds=\frac12a_1(\mathbf{u},\mathbf{v},\mathbf{w})-\frac12a_1(\mathbf{u},\mathbf{w},\mathbf{v}),\ \mathbf{u},\mathbf{v},\mathbf{w}\in \mathbf{X}\label{eq2.2}
\end{align}
by
$$ N:=\sup_{\bf{u,v,w}\in \bf{X}}\frac{b(\mathbf{u},\mathbf{v},\mathbf{w})}{|\mathbf{u}|_1 |\mathbf{v}|_1 |\mathbf{w}|_1}.$$
\begin{lemma} \label{L2.0}\cite{girault,heywood1,kos,te1} There exists a unique
solution $(\mathbf{v},q)\in (\mathbf{X},M)$
 to the steady Stokes problem
$$
-\Delta \mathbf{v}+\nabla q=\mathbf{g}, \quad \mbox{div}\;\mathbf{v}=0\quad\mbox{ in } \Omega,~~~
\mathbf{v}|_{\partial\Omega}=0,
$$
for any prescribed $\mathbf{g}\in \mathbf{Y}$ and
$$
\|\mathbf{v}\|_{2}+\|q\|_{1}\leq c\|\mathbf{g}\|_{0},
$$
where $c>0$ is a generic constant
depending on $\Omega$.
\end{lemma}

\begin{lemma} \label{L2.1} \cite{sobolevskii,McLean1993} For any $\alpha$, $t^\ast>0$ and $\phi\in L^2(0,t^\ast)$, the following property holds
 \begin{eqnarray}
\int_0^{t^\ast}\int_0^t(t-s)^{-\beta}e^{-\alpha(t-s)}\phi(s)ds\phi(t)dt\geq 0.\nonumber
\end{eqnarray}
\end{lemma}
We then introduce the Gamma function $\Gamma(\cdot)$ defined by
\begin{equation*}
\Gamma(z):=\int_0^\infty s^{z-1}e^{-s}ds,~~z>0
\end{equation*}
that will be frequently used in the future. Note that by a simple transformation, the above equation implies
\begin{equation}\label{gambnd}
\int_0^\infty s^{z-1}e^{-\upsilon s}ds=\frac{1}{\upsilon^z}\int_0^\infty s^{z-1}e^{-s}ds=\frac{\Gamma(z)}{\upsilon^z},~~z,\upsilon>0.
\end{equation}
We finally introduce the Gronwall inequality to support the analysis.
\begin{lemma} \cite{CanEwi}\label{L4.2}  If $g,h,y,G$ are nonnegative locally integrable functions on the time interval $[0,\infty)$ such that for all $t\geq 0$ and for some $C\geq 0$
\begin{eqnarray}
y(t)+G(t)\leq C+\int_{0}^th(s)ds+\int_{0}^tg(s)y(s)ds,\nonumber
\end{eqnarray}
then for $t\geq 0$
\begin{eqnarray}
y(t)+G(t)\leq \bigg(C+\int_{0}^th(s)ds\bigg)\exp\bigg(\int_{0}^tg(s)ds\bigg).\nonumber
\end{eqnarray}
\end{lemma}

\subsection{Regularity of solutions to the steady state}
By (\ref{eq2.01}) and (\ref{eq2.2}), the variational formulation of \eqref{eq1.1} reads: find $(\mathbf{\bar{u}},\bar{p})\in (\mathbf{X},M)$
such that for any $(\mathbf{v},q)\in (\mathbf{X},M)$
\begin{eqnarray}
\bigg(\mu+\frac{\rho\Gamma(1-\beta)}{\delta^{1-\beta}}\bigg) a(\bar{\mathbf{u}},\mathbf{v})-d(\mathbf{v},\bar{p})+d(\bar{\mathbf{u}},q)+a_1(\bar{\mathbf{u}},\bar{\mathbf{u}},\mathbf{v})=(\bar{\mathbf{f}},\mathbf{v}).\label{eq2.1}
  \end{eqnarray}

Then we follow \cite{He2013} to assume that the solution $\bar{\mathbf{u}}$ of (\ref{eq2.1}) satisfies
\begin{eqnarray}\label{A2}
\mu a(\mathbf{v},\mathbf{v})+b(\mathbf{v},\bar{\mathbf{u}},\mathbf{v})\geq \mu_0|\mathbf{v}|_1^2,\ \ \ \forall \ \mathbf{v}\in \mathbf{X},
\end{eqnarray}
for some $0<\mu_0<\mu$. We show in the next theorem that this assumption could ensure the uniqueness of solutions to model (\ref{eq2.1}).

\begin{theorem} \label{T3.1} Under the assumption (\ref{A2}), problem (\ref{eq2.1}) admits a unique solution pair $(\bar{\mathbf{u}},\bar{p})\in \mathbf{X}\times M$ which satisfies
\begin{eqnarray}
|\bar{\mathbf{u}}|_1\leq \bigg(\mu+\frac{\rho\Gamma(1-\beta)}{\delta^{1-\beta}}\bigg)^{-1}\|\bar{\mathbf{f}}\|_{-1},\ \ \ \|\bar{\mathbf{f}}\|_{-1}:=\sup_{\mathbf{v}\in \mathbf{X}}\frac{(\bar{\mathbf{f}},\mathbf{v})}{|\mathbf{v}|_1}.\label{eq2.7}
\end{eqnarray}
Moreover, if $\bar{\mathbf{f}}\in \mathbf{Y}$, then the solution pair $(\bar{\mathbf{u}},\bar{p})\in \mathbf{D}(A)\times (H^1(\Omega)\cap M)$ satisfies
\begin{eqnarray}
&&\|A\bar{\mathbf{u}}\|_0+\bigg(\mu+\frac{\rho\Gamma(1-\beta)}{\delta^{1-\beta}}\bigg)^{-1}\|\bar{p}\|_1\nonumber\\
&&\qquad\leq c\bigg(\mu+\frac{\rho\Gamma(1-\beta)}{\delta^{1-\beta}}\bigg)^{-1}\|\bar{\mathbf{f}}\|_{0}\bigg(1+\bigg[\mu+\frac{\rho\Gamma(1-\beta)}{\delta^{1-\beta}}\bigg]^{-4}\|\bar{\mathbf{f}}\|_{0}^2\bigg).\label{eq2.8}
\end{eqnarray}
\end{theorem}
\textbf{Proof.} By Lemma \ref{L2.0}, the proof could be carried out following those of Theorems 2.1, 2.2 and 2.4 in \cite{He2013} and is thus omitted.
$\blacksquare$

\section{Analysis of viscoelastic flow}\label{sec3}
We investigate the well-posedness and regularity of the solutions $\mathbf{u}$ to the viscoelastic fluid (\ref{eq3.1}) for $\mathbf{u}_0\in \mathbf{V}$. In the rest of the paper, we may abbreviate a space-time dependent function $g(\bf{x},t)$ as $g(t)$ for simplicity and use $\kappa$ to denote a generic positive constant that may assume different values at different occurrences.
We introduce
\begin{align*}
J(t;\mathbf{v},\mathbf{w}):&=\bigg(
\rho\int_0^t(t-s)^{-\beta}e^{-\delta(t-s)}\nabla \mathbf{v}(s)ds,\nabla\mathbf{w}(t)\bigg),~~\forall \ \mathbf{v}(\cdot,t), \mathbf{w}(\cdot,t)\in\mathbf{X}
\end{align*}
such that the variational formulation of (\ref{eq3.1}) could be formulated as
\begin{align}
&(\mathbf{u}_t,\mathbf{v})+\mu a(\mathbf{u},\mathbf{v})+a_1(\mathbf{u},\mathbf{u},\mathbf{v})-d(\mathbf{v},p)+d(\mathbf{u},q)+J(t;\mathbf{u},\mathbf{v})=(\mathbf{f},\mathbf{v}),\label{eq3.7f}\\
&\qquad\qquad\qquad\mathbf{u}_0=\mathbf{u}(t_0)\in \mathbf{V},~~\forall \ (\mathbf{v},q)\in \mathbf{X}\times M.\label{eq3.8f}\nonumber
\end{align}

\begin{lemma} \label{L4.0} For $r=0,1$, the following estimate holds
\begin{eqnarray}
\bigg|J(t;\mathbf{v},A^r\mathbf{w})\bigg|\leq\frac{\mu}{\varepsilon}\|A^{\frac{r+1}{2}}\mathbf{w}\|_0^2+\frac{\varepsilon\rho^2}{4\mu}\frac{\Gamma(1-\beta)}{\delta^{1-\beta}}\int_{0}^{t}(t-s)^{-\beta}e^{-\delta(t-s)}\|A^{\frac{r+1}{2}}\mathbf{v}\|_0^2ds\nonumber
\end{eqnarray}
for any $\mathbf{v}, \mathbf{w}\in \mathbf{V}$ where $\varepsilon>0$ is a generic constant.
\end{lemma}
\textbf{Proof.} By Cauchy-Schwarz inequality, Minkowski's integral inequality and noting $(A\mathbf{v},\mathbf{w})=(\nabla\mathbf{v},\nabla\mathbf{w})$ we have
\begin{eqnarray*}
&&\hspace{-0.1in}\bigg|J(t;\mathbf{v},A^r\mathbf{w})\bigg|\nonumber\\
&&=\bigg|\bigg(
\rho\int_0^t(t-s)^{-\beta}e^{-\delta(t-s)}A \mathbf{v}(s)ds,A^r\mathbf{w}(t)\bigg)\bigg|\nonumber\\
&&\leq\bigg\|\rho\int_0^t(t-s)^{-\beta}e^{-\delta(t-s)}A ^{\frac{r+1}{2}}\mathbf{v}(s)ds\bigg\|_0\|A^{\frac{r+1}{2}}\mathbf{w}(t)\|_0\nonumber\\
&&\leq\rho\int_0^t(t-s)^{-\beta}e^{-\delta(t-s)}\|A^{\frac{r+1}{2}} \mathbf{v}(s)\|_0ds\|A^{\frac{r+1}{2}}\mathbf{w}(t)\|_0\nonumber\\
&&\leq\frac{\mu}{\varepsilon}\|A^{\frac{r+1}{2}}\mathbf{w}(t)\|_0^2+\frac{\varepsilon}{4\mu}\bigg(\rho\int_0^t(t-s)^{-\beta}e^{-\delta(t-s)}\|A^{\frac{r+1}{2}} \mathbf{v}(s)\|_0ds\bigg)^2\nonumber\\
&&=:\frac{\mu}{\varepsilon}\|A^{\frac{r+1}{2}}\mathbf{w}(t)\|_0^2+I_0.
\end{eqnarray*}
Using Holder inequality and (\ref{gambnd}) we bound $I_0$ as
\begin{align*}
I_0&\leq\frac{\varepsilon\rho^2}{4\mu}\int_0^t(t-s)^{-\beta}e^{-\delta(t-s)}ds\int_0^t(t-s)^{-\beta}e^{-\delta(t-s)}\|A^{\frac{r+1}{2}} \mathbf{v}(s)\|_0^2ds\nonumber\\
&\leq\frac{\varepsilon\rho^2}{4\mu}\frac{\Gamma(1-\beta)}{\delta^{1-\beta}}\int_0^t(t-s)^{-\beta}e^{-\delta(t-s)}\|A^{\frac{r+1}{2}} \mathbf{v}(s)\|_0^2ds,
\end{align*}
Combining the above two equations we complete the proof.
 $\blacksquare$

\begin{theorem} Suppose $\mathbf{u}_0\in \mathbf{H}$ and $\mathbf{f}\in L^\infty(0,\infty;\mathbf{H})$. Then the viscoelastic fluid (\ref{eq3.1}) admits a unique solution $\mathbf{u}\in L^\infty(0,\infty;\mathbf{H})\cap L^2(0,t;\mathbf{V})$ for any $t\geq 0$.
\end{theorem}
\textbf{Proof.} According to the positivity of the kernel in Lemma \ref{L2.1}, the proof could be carried out following Theorem 3.1, Lemma 3.1, and Lemma 3.2 in \cite{CanEwi} and is thus omitted.
$\blacksquare$

\begin{theorem} \label{T4.1}Suppose $\mathbf{u}_0\in \mathbf{V}$, $\mathbf{f}\in L^\infty(0,\infty;\mathbf{H})\cap L^2(0,\infty;L^2(\Omega)^2)$ and satisfies
\begin{equation}\label{assf}
\int_{0}^\infty e^{2\alpha s}\|\mathbf{f}(s)\|_0^2ds<\infty.
\end{equation}
If $(\mathbf{u},p)$ satisfies (\ref{eq3.1}), then the following properties hold
\begin{eqnarray}
\|\mathbf{u}\|_0^2+\mu e^{-2\alpha t}\int_{0}^te^{2\alpha s}|\mathbf{u}|_1^2ds\leq \kappa e^{-2\alpha t},\ \ \ \forall \ t\geq 0,\label{eq3.4f}\\
|\mathbf{u}|_1^2+\mu e^{-2\alpha t}\int_{0}^te^{2\alpha s}\|A\mathbf{u}\|_0^2ds\leq\kappa e^{-2\alpha t},\ \ \ \forall \ t\geq 0,\label{eq3.5f}
\end{eqnarray}
where $\kappa>0$ depends on the data
and $0<\alpha<\frac{1}{2}\min\{\delta,\frac{\mu_0\gamma_0}{2}\}$.
\end{theorem}
\textbf{Proof.}
Taking $(\mathbf{v},q)=e^{2\alpha t}(\mathbf{u},p)$ in (\ref{eq3.7f}), we have
\begin{eqnarray}
&&\frac{e^{2\alpha t}}{2}\frac{d}{dt}\|\mathbf{u}\|_0^2+e^{2\alpha t}\mu a(\mathbf{u},\mathbf{u})+e^{2\alpha t}a_1(\mathbf{u},\mathbf{u},\mathbf{u})\nonumber\\
&&\qquad\qquad\qquad\qquad\qquad+e^{2\alpha t}J(t;\mathbf{u},\mathbf{u})=e^{2\alpha t}(\mathbf{f}(t),\mathbf{u}).\label{eq3.9f}
\end{eqnarray}
Combining (\ref{eq3.9f}) with (\ref{eq2.4}), and noting $0<\alpha\leq \frac{1}{2}\min\{\delta,\frac{\mu_0\gamma_0}{2}\}$, we obtain
\begin{eqnarray}
e^{2\alpha t}\frac{d}{dt}\|\mathbf{u}\|_0^2+2\alpha e^{2\alpha t}\|\mathbf{u}\|_0^2+\frac{3}{2}\mu e^{2\alpha t}|\mathbf{u}|_1^2+2e^{2\alpha t}J(t;\mathbf{u},\mathbf{u})\leq 2e^{2\alpha t}(\mathbf{f}(t),\mathbf{u}).\label{eq3.10f}
\end{eqnarray}
We invoke
\begin{eqnarray}
2|e^{2\alpha t}(\mathbf{f}(t),\mathbf{u})|\leq\frac{1}{2}\mu e^{2\alpha t}|\mathbf{u}|_1^2+2\mu^{-1}e^{2\alpha t}\|\mathbf{f}(t)\|_{-1}^2\label{eq3.11ad0f}\nonumber
\end{eqnarray}
in (\ref{eq3.10f}) to obtain
\begin{eqnarray}
\frac{d}{dt}(e^{2\alpha t}\|\mathbf{u}\|_0^2)+\mu e^{2\alpha t}|\mathbf{u}|_1^2+2e^{2\alpha t}J(t;\mathbf{u},\mathbf{u})\leq 2\mu^{-1}e^{2\alpha t}\|\mathbf{f}(t)\|_{-1}^2,\ \forall \ t\geq 0.\label{eq3.11ff}
\end{eqnarray}
Integrating (\ref{eq3.11ff}) from $0$ to $t$, using Lemma \ref{L2.1} and the assumption (\ref{assf}), and multiplying the resulting equation by $e^{-2\alpha t}$ on both sides we obtain
\begin{align}
\|\mathbf{u}\|_0^2+\mu e^{-2\alpha t}\int_{0}^te^{2\alpha s}|\mathbf{u}|_1^2ds&\leq e^{-2\alpha t}\|\mathbf{u}(0)\|_0^2+2\mu^{-1}e^{-2\alpha t}\int_{0}^te^{2\alpha s}\|\mathbf{f}(s)\|_0^2ds\nonumber\\
&\leq e^{-2\alpha t},\ \forall \ t\geq 0.\label{eq3.11f}
\end{align}
Then we compete the proof of (\ref{eq3.4f}).

Next, taking $(\mathbf{v},q)=e^{2\alpha t}(A\mathbf{u},0)$ in (\ref{eq3.7f}) leads to
\begin{eqnarray}
&&\frac{e^{2\alpha t}}{2}\frac{d}{dt}|\mathbf{u}|_1^2+\mu e^{2\alpha t}\|A\mathbf{u}\|_0^2+e^{2\alpha t}a_1(\mathbf{u},\mathbf{u},A\mathbf{u})\nonumber\\
&&\qquad\qquad\qquad\qquad\qquad+e^{2\alpha t}J(t;\mathbf{u},A\mathbf{u})=e^{2\alpha t}(\mathbf{f}(t),A\mathbf{u}).\label{eq3.12f}
\end{eqnarray}
By (\ref{eq2.04}) and the Young's inequality we have
\begin{align}
e^{2\alpha t}|a_1(\mathbf{u},\mathbf{u},A\mathbf{u})|&\leq c_0e^{2\alpha t}\|\mathbf{u}\|_0^\frac{1}{2}|\mathbf{u}|_1\|A\mathbf{u}\|_0^\frac{1}{2}\|A\mathbf{u}\|_0\nonumber\\
&\leq\frac{\mu}{8}e^{2\alpha t}\|A\mathbf{u}\|_0^2+8^3c_0^4\mu^{-3}e^{2\alpha t}\|\mathbf{u}\|_0^2|\mathbf{u}|_1^4\nonumber\\
&\leq\frac{\mu}{8}e^{2\alpha t}\|A\mathbf{u}\|_0^2+8^3c_0^4\mu^{-3}e^{-4\alpha t}\|e^{\alpha t}\mathbf{u}\|_0^2|e^{\alpha t}\mathbf{u}|_1^4,\nonumber\\
|e^{2\alpha t}(\mathbf{f}(t),A\mathbf{u})|&\leq\frac{\mu}{8} e^{2\alpha t}\|A\mathbf{u}\|_0^2+2\mu^{-1}e^{2\alpha t}\|\mathbf{f}(t)\|_0^2.\nonumber
\end{align}
Combining these inequalities with (\ref{eq3.12f}) and using $0<\alpha\leq \min\{\frac{\delta}{2},\frac{\mu_0\gamma_0}{4}\}$ and (\ref{eq2.4}) we get
\begin{eqnarray}
&&\frac{d}{dt}(e^{2\alpha t}|\mathbf{u}|_1^2)+\frac{5\mu}{4} e^{2\alpha t}\|A\mathbf{u}\|_0^2+2e^{2\alpha t}J(t;\mathbf{u},A\mathbf{u})\nonumber\\
&&\qquad\qquad\qquad\qquad\leq \kappa e^{-4\alpha t}|e^{\alpha t}\mathbf{u}|_1^4+2\mu^{-1}e^{2\alpha t}\|\mathbf{f}(t)\|_0^2.\label{eq3.13f}
\end{eqnarray}
Integrating (\ref{eq3.13f}) from $0$ to $t$, and noting $(A\mathbf{u},\mathbf{u})=(\nabla\mathbf{u},\nabla\mathbf{u})$, then using Lemma \ref{L2.1} yields
\begin{align}
e^{2\alpha t}|\mathbf{u}|_1^2+\mu\int_{0}^te^{2\alpha s}\|A\mathbf{u}\|_0^2ds&\leq |\mathbf{u}(0)|_1^2+\kappa \int_{0}^te^{-4\alpha s}|e^{\alpha s}\mathbf{u}|_1^2|e^{\alpha s}\mathbf{u}|_1^2ds\nonumber\\
&\qquad+2\mu^{-1}\int_{0}^te^{2\alpha s}\|\mathbf{f}(s)\|_0^2ds,\ \ \ \forall t\geq 0.\label{eq3.14f}
\end{align}
Invoking Lemma \ref{L4.2} and (\ref{eq3.11f}) in (\ref{eq3.14f}) and multiplying the resulting equation by $e^{-2\alpha t}$ on both sides we have
\begin{align}
&|\mathbf{u}|_1^2+\mu e^{-2\alpha t}\int_{0}^te^{2\alpha s}\|A\mathbf{u}\|_0^2ds\nonumber\\
&\quad\leq e^{-2\alpha t}\bigg(|\mathbf{u}(0)|_1^2+2\mu^{-1}\int_{0}^te^{2\alpha s}\|\mathbf{f}(s)\|_0^2ds\bigg)\exp\bigg(\kappa\int_{0}^te^{-4\alpha s}|e^{\alpha s}\mathbf{u}|_1^2ds\bigg)\nonumber\\
&\quad\leq e^{-2\alpha t}\bigg(|\mathbf{u}(0)|_1^2+2\mu^{-1}\int_{0}^te^{2\alpha s}\|\mathbf{f}(s)\|_0^2ds\bigg)\exp\bigg(\kappa\int_{0}^te^{-2\alpha s}|\mathbf{u}|_1^2ds\bigg)\nonumber\\
&\quad\leq \kappa e^{-2\alpha t}\bigg(|\mathbf{u}(0)|_1^2+2\mu^{-1}\int_{0}^te^{2\alpha s}\|\mathbf{f}(s)\|_0^2ds\bigg).\label{eq3.14adf}
\end{align}
Combining (\ref{eq3.14adf}) with (\ref{assf}) we obtain (\ref{eq3.5f}) and thus complete the proof of this theorem.
 $\blacksquare$

\begin{theorem} \label{T4.2} Suppose $\mathbf{u}_0\in \mathbf{D}(A)$, $\mathbf{f}\in L^\infty(0,\infty;\mathbf{H})\cap L^2(0,\infty;L^2(\Omega)^2)$, $\mathbf{f}_t\in L^2(0,\infty;H^{-1}(\Omega)^2)$ and
\begin{eqnarray}
\int_{0}^\infty e^{2\alpha s}\|\mathbf{f}(s)\|_0^2ds+\int_{0}^\infty e^{2\alpha s}\|\mathbf{f}_s(s)\|_{-1}^2ds<\infty.\label{eq3.5aadfg}
\end{eqnarray}
If $(\mathbf{u},p)$ satisfies (\ref{eq3.1}), then
\begin{align}
\|\mathbf{u}_t\|_0^2+\mu e^{-2\alpha t}\int_{0}^te^{2\alpha s}|\mathbf{u}_s|_1^2ds&\leq \kappa e^{-2\alpha t},\label{eq3.25fad}~~\|A\mathbf{u}\|_0^2\leq \kappa,\ \ \ \forall \ t\geq 0,
\end{align}
where $\kappa>0$ depends on the data
, $0<\alpha<\frac{1}{2}\min\{\delta,\frac{\mu_0\gamma_0}{2}\}$, and $\|\mathbf{f}_t(t)\|_{-1}:=\sup_{0\neq \mathbf{v}\in \mathbf{X}}\frac{(\mathbf{f}_t(t),\mathbf{v})}{|\mathbf{v}|_1}$.

In particular, if further $\|f(t)\|_0\leq c^* e^{-\alpha t}$ for $t\geq 0$ and for some constant $c^*\geq 0$, then the following exponential decay holds
\begin{equation}\label{MH}
\|\mathbf{u}\|_2^2\leq c\|A\mathbf{u}\|_0^2\leq \kappa e^{-2\alpha t},\ \ \ \forall \ t\geq 0
\end{equation}
where $c$ is given in (\ref{eq2.3}).
\end{theorem}
\textbf{Proof.}
Differentiating (\ref{eq3.7f}) with respect to time and noting
\begin{align*}
J_t(t;\mathbf{u},\mathbf{v})&=(\rho t^{-\beta}e^{-\delta t}\nabla \mathbf{u}(0),\nabla\mathbf{v})
+\bigg(\rho\int_0^t(t-s)^{-\beta}e^{-\delta(t-s)}\nabla \mathbf{u}_t(s)ds,\nabla\mathbf{v}\bigg)\nonumber\\
&=\rho t^{-\beta}e^{-\delta t}a(\mathbf{u}(0),\mathbf{v})+J(t;\mathbf{u}_t,\mathbf{v})
\end{align*}
we have
\begin{align}
&(\mathbf{u}_{tt},\mathbf{v})+\mu a(\mathbf{u}_t,\mathbf{v})+a_1(\mathbf{u}_t,\mathbf{u},\mathbf{v})+a_1(\mathbf{u},\mathbf{u}_t,\mathbf{v})-d(\mathbf{v},p_t)\nonumber\\
&\qquad+d(\mathbf{u}_t,q)+J(t;\mathbf{u}_t,\mathbf{v})=-\rho t^{-\beta}e^{-\delta t}a(\mathbf{u}(0),\mathbf{v})+(\mathbf{f}_t(t),\mathbf{v}).\label{eq3.21f}
\end{align}
Taking $(\mathbf{v},q)=e^{2\alpha t}(\mathbf{u}_t,p_t)$ in (\ref{eq3.21f}) we get
\begin{align}
&\frac{e^{2\alpha t}}{2}\frac{d}{dt}\|\mathbf{u}_t\|_0^2+\mu e^{2\alpha t}|\mathbf{u}_t|_1^2+e^{2\alpha t}a_1(\mathbf{u}_t,\mathbf{u},\mathbf{u}_t)+e^{2\alpha t} J(t;\mathbf{u}_t,\mathbf{u}_t)\nonumber\\
&\qquad=-\rho t^{-\beta}e^{-(\delta-2\alpha) t}a(\mathbf{u}(0),\mathbf{u}_t)+e^{2\alpha t}(\mathbf{f}_t(t),\mathbf{u}_t).\label{eq3.22f}
\end{align}
In view of (\ref{eq2.03}) and $\mathbf{u}_0\in \mathbf{D}(A)$ we obtain
\begin{align}
\rho t^{-\beta}e^{-(\delta-2\alpha) t}|-a(\mathbf{u}(0),\mathbf{u}_t)|&\leq \rho t^{-\beta}e^{-(\delta-2\alpha) t}\|A\mathbf{u}(0)\|_0\|\mathbf{u}_t\|_0\nonumber\\
&\leq \kappa t^{-\beta}e^{-(\delta-2\alpha) t}\|\mathbf{u}_t\|_0,\nonumber\\
e^{2\alpha t}|a_1(\mathbf{u}_t,\mathbf{u},\mathbf{u}_t)|&\leq c_0e^{2\alpha t}|\mathbf{u}_t|_1\|A\mathbf{u}\|_0\|\mathbf{u}_t\|_0\nonumber\\
&\leq \frac{\mu}{8}e^{2\alpha t}|\mathbf{u}_t|_1^2+2c_0^2\mu^{-1}e^{2\alpha t}\|A\mathbf{u}\|_0^2\|\mathbf{u}_t\|_0^2,\nonumber\\
e^{2\alpha t}(\mathbf{f}_t(t),\mathbf{u}_t)&\leq\frac{\mu}{8}e^{2\alpha t}|\mathbf{u}_t|_1^2+2\mu^{-1}e^{2\alpha t}\|\mathbf{f}_t(t)\|_{-1}^2.\nonumber
\end{align}
Combining these inequalities with (\ref{eq3.22f}) and noting $0<\alpha\leq \min\{\frac{\delta}{2},\frac{\mu_0\gamma_0}{4}\}$ yield
\begin{align}
&e^{2\alpha t}\frac{d}{dt}\|\mathbf{u}_t\|_0^2+2\alpha e^{2\alpha t}\|\mathbf{u}_t\|_0^2+\mu e^{2\alpha t}|\mathbf{u}_t|_1^2+2e^{2\alpha t} J(t;\mathbf{u}_t,\mathbf{u}_t)\nonumber\\
&\quad\leq 2c_0^2\mu^{-1}e^{2\alpha t}\|A\mathbf{u}\|_0^2\|\mathbf{u}_t\|_0^2+\kappa t^{-\beta}e^{-(\delta-2\alpha) t}\|\mathbf{u}_t\|_0+2\mu^{-1}e^{2\alpha t}\|\mathbf{f}_t(t)\|_{-1}^2.\label{eq3.23f}
\end{align}
We integrate (\ref{eq3.23f}) from $0$ to $t$, use $\|\mathbf{u}_t(0)\|_0^2\leq c\|A\mathbf{u}(0)\|_0^2$ derived from model (\ref{eq3.1}), and $(A\mathbf{u_t},\mathbf{u_t})=(\nabla\mathbf{u_t},\nabla\mathbf{u_t})$ as well as Lemma \ref{L2.1} and (\ref{eq3.5aadfg}) to obtain
\begin{align}
&e^{2\alpha t}\|\mathbf{u}_t\|_0^2+\mu \int_{0}^te^{2\alpha s}|\mathbf{u}_s|_1^2ds\nonumber\\
&\leq \|\mathbf{u}_t(0)\|_0^2+2c_0^2\mu^{-1}\int_{0}^te^{2\alpha s}\|A\mathbf{u}\|_0^2\|\mathbf{u}_s\|_0^2ds\nonumber\\
&\qquad+\kappa \int_{0}^ts^{-\beta}e^{-(\delta-2\alpha) s}\|\mathbf{u}_s\|_0ds+2\mu^{-1}\int_{0}^te^{2\alpha s}\|\mathbf{f}_s(s)\|_{-1}^2ds\nonumber\\
&\leq\kappa+2c_0^2\mu^{-1}\int_{0}^t\|A\mathbf{u}\|_0^2e^{2\alpha s}\|\mathbf{u}_s\|_0^2ds\nonumber\\
&\qquad+\kappa\int_{0}^ts^{-\beta}e^{-(\delta-\alpha) s}e^{2\alpha s}\|\mathbf{u}_s\|_0^2ds.\label{eq3.24f}
\end{align}
Using Lemma \ref{L4.2} and (\ref{eq3.5f}) for (\ref{eq3.24f}) we have
\begin{align*}
&e^{2\alpha t}\|\mathbf{u}_t\|_0^2+\mu \int_{0}^te^{2\alpha s}|\mathbf{u}_s|_1^2ds\nonumber\\
&\qquad\leq\kappa\exp\bigg(\kappa\int_{0}^t\|A\mathbf{u}\|_0^2ds+\kappa\int_{0}^ts^{-\beta}e^{-(\delta-\alpha) s}ds\bigg)\leq \kappa,
\end{align*}
which proves the first equation of (\ref{eq3.25fad}).

Then we take $(\mathbf{v},q)=(A\mathbf{u},0)$ in (\ref{eq3.7f}) to get
\begin{align}
(\mathbf{u}_t,A\mathbf{u})+\mu \|A\mathbf{u}\|_0^2+a_1(\mathbf{u},\mathbf{u},A\mathbf{u})+J(t;\mathbf{u},A\mathbf{u})=(\mathbf{f}(t),A\mathbf{u}).\label{eq3.12adf}
\end{align}
Due to (\ref{eq2.03}), we have
\begin{align}
|a_1(\mathbf{u},\mathbf{u},A\mathbf{u})|&\leq c_0\|\mathbf{u}\|_0^\frac{1}{2}|\mathbf{u}|_1^\frac{1}{2}|\mathbf{u}|_1^\frac{1}{2}\|A\mathbf{u}\|_0^\frac{3}{2}\nonumber\\
&\leq\frac{\mu}{6}\|A\mathbf{u}\|_0^2+6^3 c_0^4\mu^{-3}\|\mathbf{u}\|_0^2|\mathbf{u}|_1^2|\mathbf{u}|_1^2,\nonumber\\
|(\mathbf{f}(t),A\mathbf{u})|&\leq\frac{\mu}{6}\|A\mathbf{u}\|_0^2+\frac{3}{2}\mu^{-1}\|\mathbf{f}(t)\|_0^2.\nonumber
\end{align}
Combining these inequalities with (\ref{eq3.12adf}) and using (\ref{eq3.4f}) and (\ref{eq3.5f}) yield
\begin{align}
&2(\mathbf{u}_t,A\mathbf{u})+\frac{4}{3}\mu \|A\mathbf{u}\|_0^2\leq \kappa |\mathbf{u}|_1^2+\frac{3}{2}\mu^{-1}\|\mathbf{f}(t)\|_0^2-2J(t;\mathbf{u},A\mathbf{u}).\label{eq3.28}
\end{align}
Using Lemma \ref{L4.0}, we get
\begin{eqnarray}
\bigg|J(t;\mathbf{u},A\mathbf{u})\bigg|\leq\frac{\mu}{12}\|A\mathbf{u}\|_0^2+\frac{3\rho^2}{\mu}\frac{\Gamma(1-\beta)}{\delta^{1-\beta}}\int_{0}^{t}(t-s)^{-\beta}e^{-\delta(t-s)}\|A\mathbf{u}\|_0^2ds.\label{eq7.0}
\end{eqnarray}
We incorporate
\begin{eqnarray}
&&2|(\mathbf{u}_t,A\mathbf{u})|\leq\frac{\mu}{6}\|A\mathbf{u}\|_0^2+6\mu^{-1}\|\mathbf{u}_t\|_0^2\nonumber
\end{eqnarray}
with (\ref{eq3.28}) and use (\ref{eq7.0}) to obtain
\begin{eqnarray}
&&\mu \|A\mathbf{u}\|_0^2\leq \kappa |\mathbf{u}|_1^2+6\mu^{-1}\|\mathbf{u}_t\|_0^2+\frac{3}{2}\mu^{-1}\|\mathbf{f}(t)\|_0^2\nonumber\\
&&\qquad\qquad+\frac{6\rho^2}{\mu}\frac{\Gamma(1-\beta)}{\delta^{1-\beta}}\int_{0}^{t}(t-s)^{-\beta}e^{-\delta(t-s)}\|A\mathbf{u}\|_0^2ds.\label{eq7.1}
\end{eqnarray}
Applying (\ref{eq3.25fad}), (\ref{eq3.5f}) and (\ref{eq2.3}) in (\ref{eq7.1}) we get
\begin{eqnarray}
\mu \|A\mathbf{u}\|_0^2\leq \kappa (e^{-2\alpha t}+\|\mathbf{f}(t)\|_0^2) +\frac{6\rho^2}{\mu}\frac{\Gamma(1-\beta)}{\delta^{1-\beta}}\int_{0}^{t}(t-s)^{-\beta}e^{-\delta(t-s)}\|A\mathbf{u}\|_0^2ds.\nonumber
\end{eqnarray}
By Lemma \ref{L4.2}, we have
\begin{eqnarray}\label{zx}
\mu \|A\mathbf{u}\|_0^2\leq \kappa (e^{-2\alpha t}+\|\mathbf{f}(t)\|_0^2) \exp\bigg(\frac{6\rho^2}{\mu}\frac{\Gamma(1-\beta)}{\delta^{1-\beta}}\int_{0}^{t}(t-s)^{-\beta}e^{-\delta(t-s)}ds\bigg),
\end{eqnarray}
which, together with (\ref{gambnd}) and $f\in L^\infty(0,\infty;\mathbf{H})$, leads to the second estimate of (\ref{eq3.25fad}). We finally use (\ref{zx}) and $\|f(t)\|_0\leq c^* e^{-\alpha t}$ to obtain (\ref{MH}) and thus complete the proof.
 $\blacksquare$

\section{Exponential convergence: basic results}\label{sec4}
We prove basic asymptotic results on the solutions between (\ref{eq3.1}) and (\ref{eq1.1}). Let $\mathbf{z}:=\mathbf{u}-\bar{\mathbf{u}}$, $\eta:=p-\bar{p}$ and
\begin{eqnarray}
J_1(t;\bar{\mathbf{u}},\mathbf{v}):=\frac{\rho}{\delta^{1-\beta}}\int_{\delta t}^\infty s^{-\beta}e^{-s}ds \,a(\bar{\mathbf{u}},\mathbf{v}(t)),~~\forall \ \mathbf{v}(\cdot,t)\in\mathbf{X}.\nonumber
\end{eqnarray}


\begin{lemma} \label{L5.1} The following estimates hold for $r=0,1$
\begin{eqnarray}
e^{2\alpha t}|J_1(t;\bar{\mathbf{u}},A^r\mathbf{v})|\leq c_2 e^{-(\delta-2\alpha) t}\|A^\frac{r+1}{2}\bar{\mathbf{u}}\|_{0}^2+\frac{\mu}{4\varepsilon}e^{2\alpha t}\|A^\frac{r+1}{2}\mathbf{v}\|_{0}^2,\ \ \forall \ \mathbf{v}\in \mathbf{V},\nonumber
\end{eqnarray}
where $\varepsilon>0$ is a general constant, $c_2>0$ depends only on the data and $\bar{\mathbf{u}}\in \mathbf{V}$.
\end{lemma}
\textbf{Proof.} Since
\begin{align}
\frac{\rho}{\delta^{1-\beta}}\int_{\delta t}^\infty s^{-\beta}e^{-s}ds&\leq \frac{\rho}{\delta^{1-\beta}}e^{-\frac{\delta t}{2}}\int_{\delta t}^\infty s^{-\beta}e^{-\frac{s}{2}}ds\leq c_1e^{-\frac{\delta t}{2}},\label{eq3.11ad1}
\end{align}
where $c_1>0$ is a constant depends only on the data $\rho,\delta,\beta$, noting $(A\bar{\mathbf{u}},\mathbf{v})=(\nabla\bar{\mathbf{u}},\nabla\mathbf{v})$, we have for $\mathbf{v}\in \mathbf{V}$
\begin{align}
&|e^{2\alpha t}J_1(t;\bar{\mathbf{u}},A^r\mathbf{v})|\nonumber\\
&\qquad=e^{2\alpha t}\frac{\rho}{\delta^{1-\beta}}\int_{\delta t}^\infty s^{-\beta}e^{-s}ds|(A\bar{\mathbf{u}},A^r\mathbf{v})|\nonumber\\
&\qquad\leq \varepsilon\mu_0^{-1}e^{2\alpha t}\bigg(\frac{\rho}{\delta^{1-\beta}}\int_{\delta t}^\infty s^{-\beta}e^{-s}ds\bigg)^2\|A^\frac{r+1}{2}\bar{\mathbf{u}}\|_{0}^2+\frac{\mu}{4\varepsilon}e^{2\alpha t}\|A^\frac{r+1}{2}\mathbf{v}\|_{0}^2\nonumber\\
&\qquad\leq c_2 e^{-(\delta-2\alpha) t}\|A^\frac{r+1}{2}\bar{\mathbf{u}}\|_{0}^2+\frac{\mu}{4\varepsilon}e^{2\alpha t}\|A^\frac{r+1}{2}\mathbf{v}\|_{0}^2,\nonumber
\end{align}
which completes the proof.
 $\blacksquare$

\begin{theorem} \label{T5.1} Suppose the assumptions of the Theorem \ref{T4.2} hold, $\bar{\mathbf{f}}\in L^2(\Omega)^2$ and
\begin{eqnarray}
\int_{0}^\infty e^{2\alpha s}\|\mathbf{f}(s)-\bar{\mathbf{f}}\|_0^2ds<\infty.\label{eq3.5aad}
\end{eqnarray}
Then the solution $\mathbf{u}$ of (\ref{eq3.1}) exponentially converges to $\bar{\mathbf{u}}$ in the following sense
\begin{eqnarray}
\|\mathbf{z}\|_0^2+\mu e^{-2\alpha t}\int_{0}^te^{2\alpha s}|\mathbf{z}|_1^2ds\leq \kappa e^{-2\alpha t},\ \ \ \forall \ t\geq 0,\label{eq3.4}\\
|\mathbf{z}|_1^2+\mu e^{-2\alpha t}\int_{0}^te^{2\alpha s}\|A\mathbf{z}\|_0^2ds\leq\kappa e^{-2\alpha t},\ \ \ \forall \ t\geq 0,\label{eq3.5}
\end{eqnarray}
where $\kappa>0$ depends on the data
 and $0<\alpha<\frac{1}{2}\min\{\delta,\frac{\mu_0\gamma_0}{2}\}$.
\end{theorem}
\textbf{Proof.} By (\ref{gambnd}) the equation (\ref{eq1.1}) could be reformulated as follows
\begin{align}
&\bar{\mathbf{u}}_t-\mu\Delta \bar{\mathbf{u}}+(\bar{\mathbf{u}}\cdot\nabla)\bar{\mathbf{u}}+\nabla
\bar{p}-
\rho\int_0^t(t-s)^{-\beta}e^{-\delta(t-s)}\Delta \bar{\mathbf{u}}ds\nonumber\\
&\qquad\qquad=\bar{\mathbf{f}}+
\frac{\rho}{\delta^{1-\beta}}\int_{\delta t}^\infty s^{-\beta}e^{-s}\Delta \bar{\mathbf{u}}ds.\label{eq5.0}
\end{align}
Subtracting (\ref{eq5.0}) from (\ref{eq3.1}) we obtain that
\begin{eqnarray}
&&\mathbf{z}_t-\mu\Delta \mathbf{z}+(\mathbf{z}\cdot\nabla)\mathbf{z}+(\bar{\mathbf{u}}\cdot\nabla)\mathbf{z}+(\mathbf{z}\cdot\nabla)\bar{\mathbf{u}}-
\rho\int_0^t(t-s)^{-\beta}e^{-\delta(t-s)}\Delta \mathbf{z}(s)ds\nonumber\\
&&\qquad+\nabla
\eta=-\frac{\rho}{\delta^{1-\beta}}\int_{\delta t}^\infty s^{-\beta}e^{-s}\Delta \bar{\mathbf{u}}(s)ds+\mathbf{f}(t)-\bar{\mathbf{\mathbf{f}}},\label{eq3.6}\\
&&\mathbf{z}(t)\in \mathbf{V}, \ \eta(t)\in M,\ t\geq 0, \ \  \mathbf{z}(0)=\mathbf{u}_0-\bar{\mathbf{u}}\in \mathbf{V},\label{eq3.8a}
\end{eqnarray}
for all $(\mathbf{x},t)\in\Omega\times (0,\infty)$. Then the variational formulation of (\ref{eq3.6}) are
\begin{eqnarray}
&&(\mathbf{z}_t,\mathbf{v})+\mu a(\mathbf{z},\mathbf{v})+a_1(\mathbf{z},\mathbf{z},\mathbf{v})+a_1(\bar{\mathbf{u}},\mathbf{z},\mathbf{v})+a_1(\mathbf{z},\bar{\mathbf{u}},\mathbf{v})\nonumber\\
&&\qquad-d(\mathbf{v},\eta)+d(\mathbf{z},q)+J(t;\mathbf{z},\mathbf{v})=J_1(t;\bar{\mathbf{u}},\mathbf{v})+(\mathbf{f}(t)-\bar{\mathbf{f}},\mathbf{v}),\label{eq3.7}\\
&&\mathbf{z}(0)=\mathbf{u}_0-\bar{\mathbf{u}}\in \mathbf{V},\nonumber
\end{eqnarray}
for all $(\mathbf{v},q)\in (\mathbf{X}\times M)$. Taking $(\mathbf{v},q)=e^{2\alpha t}(\mathbf{z},\eta)$ in (\ref{eq3.7}), using (\ref{eq2.01}) and noting $\mbox{div}\;\mathbf{z}=0$ from (\ref{eq3.8a}) we have
\begin{eqnarray}
&&\frac{e^{2\alpha t}}{2}\frac{d}{dt}\|\mathbf{z}\|_0^2+e^{2\alpha t}\mu a(\mathbf{z},\mathbf{z})+e^{2\alpha t}b(\mathbf{z},\bar{\mathbf{u}},\mathbf{z})+e^{2\alpha t}J(t;\mathbf{z},\mathbf{z})\nonumber\\
&&\qquad=e^{2\alpha t}J_1(t;\bar{\mathbf{u}},\mathbf{z})+e^{2\alpha t}(\mathbf{f}(t)-\bar{\mathbf{f}},\mathbf{z}).\label{eq3.9}
\end{eqnarray}
Combining the assumption (\ref{A2}), (\ref{eq2.4}) and noting $0<\alpha<\frac{1}{2}\min\{\delta,\frac{\mu_0\gamma_0}{2}\}$ we obtain from (\ref{eq3.9}) that
\begin{eqnarray}
&&e^{2\alpha t}\frac{d}{dt}\|\mathbf{z}\|_0^2+2\alpha e^{2\alpha t}\|\mathbf{z}\|_0^2+\frac{3}{2}\mu e^{2\alpha t}|\mathbf{z}|_1^2+2e^{2\alpha t}J(t;\mathbf{z},\mathbf{z})\nonumber\\
&&\qquad\leq 2e^{2\alpha t}J_1(t;\bar{\mathbf{u}},\mathbf{z})+2e^{2\alpha t}(\mathbf{f}(t)-\bar{\mathbf{f}},\mathbf{z}).\label{eq3.10}
\end{eqnarray}
By Lemma \ref{L5.1} we have
\begin{align*}
2e^{2\alpha t}|J_1(t;\bar{\mathbf{u}},\mathbf{z})|&\leq c_2 e^{-(\delta-2\alpha) t}|\bar{\mathbf{u}}|_{1}^2+\frac{\mu}{4}e^{2\alpha t}|\mathbf{z}|_{1}^2,\\
2|e^{2\alpha t}(\mathbf{f}(t)-\bar{\mathbf{f}},\mathbf{z})|&\leq\frac{\mu}{4} e^{2\alpha t}|\mathbf{z}|_1^2+4\mu^{-1}e^{2\alpha t}\|\mathbf{f}(t)-\bar{\mathbf{f}}\|_0^2,
\end{align*}
which, together with (\ref{eq3.10}), yield
\begin{align}
&\frac{d}{dt}(e^{2\alpha t}\|\mathbf{z}\|_0^2)+\mu e^{2\alpha t}|\mathbf{z}|_1^2+2e^{2\alpha t}J(t;\mathbf{z},\mathbf{z})\nonumber\\
&\qquad\leq c_2 e^{-(\delta-2\alpha) t}|\bar{\mathbf{u}}|_{1}^2+4\mu^{-1}e^{2\alpha t}\|\mathbf{f}(t)-\bar{\mathbf{f}}\|_0^2.\label{eq3.11y}
\end{align}
Integrating (\ref{eq3.11y}) from $0$ to $t$ and using Lemma \ref{L2.1}, (\ref{eq3.5aad}) we deduce
\begin{align}
&\|\mathbf{z}\|_0^2+\mu e^{-2\alpha t}\int_{0}^te^{2\alpha s}|\mathbf{z}|_1^2ds\nonumber\\
&\qquad\leq e^{-2\alpha t}\|\mathbf{z}(0)\|_0^2+c_2|\bar{\mathbf{u}}|_1^2e^{-2\alpha t}\int_{0}^te^{-(\delta-2\alpha)s}ds\nonumber\\
&\qquad\qquad+\frac{4e^{-2\alpha t}}{\mu}\int_{0}^te^{2\alpha s}\|\mathbf{f}(s)-\bar{\mathbf{f}}\|_0^2ds\nonumber\\
&\qquad\leq e^{-2\alpha t}\|\mathbf{z}(0)\|_0^2+\kappa e^{-2\alpha t}\leq \kappa e^{-2\alpha t},\ \ \ \forall \ t\geq 0,\nonumber
\end{align}
which proves (\ref{eq3.4}).

Next, taking $(\mathbf{v},q)=e^{2\alpha t}(A\mathbf{z},0)$ in (\ref{eq3.7}) we get
\begin{align}
&\mu e^{2\alpha t}\|A\mathbf{z}\|_0^2+\frac{e^{2\alpha t}}{2}\frac{d}{dt}|\mathbf{z}|_1^2+e^{2\alpha t}J(t;\mathbf{z},A\mathbf{z})\nonumber\\
&\qquad+e^{2\alpha t}a_1(\mathbf{z},\mathbf{z},A\mathbf{z})+e^{2\alpha t}a_1(\bar{\mathbf{u}},\mathbf{z},A\mathbf{z})+e^{2\alpha t}a_1(\mathbf{z},\bar{\mathbf{u}},A\mathbf{z})\nonumber\\
&\quad=e^{2\alpha t}J_1(t;\bar{\mathbf{u}},A\mathbf{z})+e^{2\alpha t}(\mathbf{f}(t)-\bar{\mathbf{f}},A\mathbf{z}).\label{eq3.12}
\end{align}
Due to (\ref{eq2.03}) and Lemma \ref{L5.1}, we have
\begin{align}
e^{2\alpha t}|a_1(\mathbf{z},\mathbf{z},A\mathbf{z})|&\leq c_0e^{2\alpha t}\|\mathbf{z}\|_0^\frac{1}{2}|\mathbf{z}|_1\|A\mathbf{z}\|_0^\frac{1}{2}\|A\mathbf{z}\|_0\nonumber\\
&\leq\frac{\mu}{12}e^{2\alpha t}\|A\mathbf{z}\|_0^2+12^3c_0^4\mu^{-3}e^{2\alpha t}\|\mathbf{z}\|_0^2|\mathbf{z}|_1^4,\nonumber\\
e^{2\alpha t}|a_1(\bar{\mathbf{u}},\mathbf{z},A\mathbf{z})|&+e^{2\alpha t}|a_1(\mathbf{z},\bar{\mathbf{u}},A\mathbf{z})|\nonumber\\
&\leq 2c_0e^{2\alpha t}|\mathbf{z}|_1\|A\bar{\mathbf{u}}\|_0\|A\mathbf{z}\|_0\nonumber\\
&\leq\frac{\mu}{12}e^{2\alpha t}\|A\mathbf{z}\|_0^2+12c_0^2\mu^{-1}e^{2\alpha t}\|A\bar{\mathbf{u}}\|_0^2|\mathbf{z}|_1^2,\nonumber\\
e^{2\alpha t}|J_1(t;\bar{\mathbf{u}},A\mathbf{z})|&\leq c_2 e^{-(\delta-2\alpha) t}\|A\bar{\mathbf{u}}\|_{0}^2+\frac{\mu}{12}e^{2\alpha t}\|A\mathbf{z}\|_{0}^2,\nonumber\\
|e^{2\alpha t}(\mathbf{f}(t)-\bar{\mathbf{f}},A\mathbf{z})|&\leq\frac{\mu}{12} e^{2\alpha t}\|A\mathbf{z}\|_0^2+3\mu^{-1}e^{2\alpha t}\|\mathbf{f}(t)-\bar{\mathbf{f}}\|_0^2.\nonumber
\end{align}
Invoking these inequalities and (\ref{eq2.4}),  (\ref{eq2.8}), (\ref{eq3.4}), $0<\alpha<\frac{1}{2}\min\{\delta,\frac{\mu_0\gamma_0}{2}\}$ as well as $|\mathbf{z}|_1\leq|\bar{\mathbf{u}}|_1+|\mathbf{u}(t)|_1\leq\kappa$ proved by (\ref{eq2.7}) and (\ref{eq3.5f}) in (\ref{eq3.12}) we get
\begin{eqnarray}
&&\frac{d}{dt}(e^{2\alpha t}|\mathbf{z}|_1^2)+\mu e^{2\alpha t}\|A\mathbf{z}\|_0^2+2e^{2\alpha t}J(t;\mathbf{z},A\mathbf{z})\nonumber\\
&&\qquad\leq \kappa e^{2\alpha t}|\mathbf{z}|_1^2+\kappa\|A\bar{\mathbf{u}}\|_{0}^2+6\mu^{-1}e^{2\alpha t}\|\mathbf{f}(t)-\bar{\mathbf{f}}\|_0^2.\label{eq3.13}
\end{eqnarray}
Integrating (\ref{eq3.13}) from $0$ to $t$ and noting $(A\mathbf{z},\mathbf{z})=(\nabla\mathbf{z},\nabla\mathbf{z})$, then using Lemma \ref{L2.1} we reach
\begin{align*}
&|\mathbf{z}|_1^2+\mu e^{-2\alpha t}\int_{0}^te^{2\alpha s}\|A\mathbf{z}\|_0^2ds\\
&\qquad\leq e^{-2\alpha t}|\mathbf{z}(0)|_1^2+\kappa e^{-2\alpha t}\int_{0}^te^{2\alpha s}|\mathbf{z}|_1^2ds+\kappa e^{-2\alpha t}\nonumber\\
&\qquad\qquad+6\mu^{-1}e^{-2\alpha t}\int_{0}^te^{2\alpha s}\|\mathbf{f}(s)-\bar{\mathbf{f}}\|_0^2ds,\ \ \ \forall t\geq 0,
\end{align*}
which, together with (\ref{eq3.5aad}) and (\ref{eq3.4}), leads to (\ref{eq3.5}).
 $\blacksquare$

\begin{theorem} Under the assumptions of the Theorem \ref{T5.1}, $p$ converges to $\bar{p}$  in an exponential rate
\begin{eqnarray}
\|p(t)-\bar{p}\|_0^2\leq\kappa e^{-2\alpha t},\ \ \ \forall \ t\geq 0,\label{eq3.17}
\end{eqnarray}
where $\kappa>0$ depends on the data
 and $0<\alpha<\frac{1}{2}\min\{\delta,\frac{\mu_0\gamma_0}{2}\}$.
\end{theorem}
\textbf{Proof.}
By $\mathbf{z}_t=\mathbf{u}_t$ and the estimate of $\mathbf{u}_t$ in Theorem \ref{T4.2} we have
\begin{eqnarray}
\|\mathbf{z}_t\|_0^2\leq \kappa e^{-2\alpha t},\ \ \forall\ t\geq0.\label{eq3.25}
\end{eqnarray}
Using (\ref{eq2.02}), (\ref{eq3.7}) and (\ref{eq5.1}), we obtain
\begin{align}
\|\eta(t)\|_0&\leq C_0^{-1}\|\mathbf{f}(t)-\bar{\mathbf{f}}\|_{0}+ C_0^{-1}\|\mathbf{z}_t\|_0+C_0^{-1}|\mathbf{z}|_1\nonumber\\
&+C_0^{-1}(|\mathbf{z}|_1+|\bar{\mathbf{u}}|_1)|\mathbf{z}|_1+C_0^{-1}\int_0^t(t-s)^{-\beta}e^{-\delta(t-s)}|\mathbf{z}|_1ds\nonumber\\
&+C_0^{-1}\frac{\rho}{\delta^{1-\beta}}\int_{\delta t}^\infty s^{-\beta}e^{-s}|\bar{\mathbf{u}}|_1ds.\label{eq3.26}
\end{align}
By (\ref{eq3.5}), (\ref{eq3.11ad1}), (\ref{eq2.7}) and $0<\alpha<\frac{1}{2}\min\{\delta,\frac{\mu_0\gamma_0}{2}\}$ we get
\begin{align}
\int_0^t(t-s)^{-\beta}e^{-\delta(t-s)}|\mathbf{z}|_1ds&\leq e^{-\alpha t}\int_0^t(t-s)^{-\beta}e^{-(\delta-\alpha)(t-s)}e^{\alpha s}|\mathbf{z}|_1ds\nonumber\\
&\leq e^{-\alpha t}\int_0^t(t-s)^{-\beta}e^{-(\delta-\alpha)(t-s)}ds\leq \kappa e^{-\alpha t},\label{eq3.26dda1}\\
\frac{\rho}{\delta^{1-\beta}}\int_{\delta t}^\infty s^{-\beta}e^{-s}|\bar{\mathbf{u}}|_1ds&\leq\kappa e^{-\frac{\delta t}{2}}\leq\kappa e^{-\alpha t}.\label{eq3.27dda1}
\end{align}
Combining (\ref{eq3.26}) with (\ref{eq3.26dda1})--(\ref{eq3.27dda1}) and using (\ref{eq3.25}), (\ref{eq3.5}), (\ref{eq3.5aad}) and (\ref{eq2.7}) we prove (\ref{eq3.17}).
 $\blacksquare$

\section{Exponential convergence in stronger norms}\label{sec5}

\begin{theorem} \label{T5.3} Under the assumptions of Theorem \ref{T5.1}, the solution $(\mathbf{u},p)$ of (\ref{eq3.1}) converges to the solution $(\bar{\mathbf{u}},\bar{p})$ of (\ref{eq1.1}) in an exponential rate
\begin{eqnarray}
\|\mathbf{u}(t)-\bar{\mathbf{u}}\|_2^2\leq\kappa e^{-2\alpha t},\ \ \ \forall \ t\geq 0,\label{eq3.15ad}\\
\|p(t)-\bar{p}\|_1^2\leq\kappa e^{-2\alpha t},\ \ \ \forall \ t\geq 0,\label{eq3.17ad}
\end{eqnarray}
where $\kappa>0$ depends on the data
 and $0<\alpha<\frac{1}{2}\min\{\delta,\frac{\mu_0\gamma_0}{2}\}$.
\end{theorem}
\textbf{Proof.}
Taking $(\mathbf{v},q)=(A\mathbf{z},0)$ in (\ref{eq3.7}) we get
\begin{align}
&(\mathbf{z}_t,A\mathbf{z})+\mu \|A\mathbf{z}\|_0^2+a_1(\mathbf{z},\mathbf{z},A\mathbf{z})+a_1(\bar{\mathbf{u}},\mathbf{z},A\mathbf{z})+a_1(\mathbf{z},\bar{\mathbf{u}},A\mathbf{z})+J(t;\mathbf{z},A\mathbf{z})\nonumber\\
&\qquad=J_1(t;\bar{\mathbf{u}},A\mathbf{z})+(\mathbf{f}(t)-\bar{\mathbf{f}},A\mathbf{z}).\label{eq3.12ad}
\end{align}
By (\ref{eq2.03}) and (\ref{eq2.04}) we derive
\begin{align}
|a_1(\mathbf{z},\mathbf{z},A\mathbf{z})|&\leq c_0\|\mathbf{z}\|_0^\frac{1}{2}|\mathbf{z}|_1^\frac{1}{2}|\mathbf{z}|_1^\frac{1}{2}\|A\mathbf{z}\|_0^\frac{3}{2},\nonumber\\
&\leq\frac{\mu}{8}\|A\mathbf{z}\|_0^2+8^3c_0^4\mu^{-3}\|\mathbf{z}\|_0^2|\mathbf{z}|_1^2|\mathbf{z}|_1^2,\nonumber\\
|a_1(\bar{\mathbf{u}},\mathbf{z},A\mathbf{z})|+|a_1(\mathbf{z},\bar{\mathbf{u}},A\mathbf{z})|&\leq 2c_0|\mathbf{z}|_1\|A\bar{\mathbf{u}}\|_0\|A\mathbf{z}\|_0,\nonumber\\
&\leq\frac{\mu}{8}\|A\mathbf{z}\|_0^2+8c_0^2\mu^{-1}\|A\bar{\mathbf{u}}\|_0^2|\mathbf{z}|_1^2,\nonumber\\
|(\mathbf{f}(t)-\bar{\mathbf{f}},A\mathbf{z})|&\leq\frac{\mu}{8}\|A\mathbf{z}\|_0^2+2\mu^{-1}\|\mathbf{f}(t)-\bar{\mathbf{f}}\|_0^2.\nonumber
\end{align}
Combining these inequalities with (\ref{eq3.12ad}) and using (\ref{eq2.8}), (\ref{eq3.4}) and (\ref{eq3.5}) we obtain
\begin{eqnarray}
&&2(\mathbf{z}_t,A\mathbf{z})+\frac{5}{4}\mu \|A\mathbf{z}\|_0^2+2J(t;\mathbf{z},A\mathbf{z})\nonumber\\
&&\qquad\leq \kappa |\mathbf{z}|_1^2+2J_1(t;\bar{\mathbf{u}},A\mathbf{z})+4\mu^{-1}\|\mathbf{f}(t)-\bar{\mathbf{f}}\|_0^2.\label{eq3.13ad}
\end{eqnarray}
Using Lemma \ref{L4.0}, we get
\begin{eqnarray}
\bigg|J(t;\mathbf{z},A\mathbf{z})\bigg|\leq\frac{\mu}{24}\|A\mathbf{z}\|_0^2+\frac{6\rho^2}{\mu}\frac{\Gamma(1-\beta)}{\delta^{1-\beta}}\int_{0}^{t}(t-s)^{-\beta}e^{-\delta(t-s)}\|A\mathbf{z}\|_0^2ds.\label{eq7.2}
\end{eqnarray}
We incorporate
\begin{align}
2|(\mathbf{z}_t,A\mathbf{z})|&\leq\frac{\mu}{12}\|A\mathbf{z}\|_0^2+12\mu^{-1}\|\mathbf{z}_t\|_0^2,\nonumber\\
2|J_1(t;\bar{\mathbf{u}},A\mathbf{z})|&=\frac{2\rho}{\delta^{1-\beta}}\int_{\delta t}^\infty s^{-\beta}e^{-s}ds|(A\bar{\mathbf{u}},A\mathbf{z})|\nonumber\\
&\leq \frac{\mu}{12}\|A\mathbf{z}\|_0^2+\frac{12}{\mu}\bigg(\frac{\rho}{\delta^{1-\beta}}\int_{\delta t}^\infty s^{-\beta}e^{-s}ds\bigg)^2\|A\bar{\mathbf{u}}\|_{0}^2.\nonumber
\end{align}
with (\ref{eq3.13ad}) and (\ref{eq7.2}) to obtain
\begin{align}
\mu \|A\mathbf{z}\|_0^2&\leq \kappa |\mathbf{z}|_1^2+12\mu^{-1}\|\mathbf{z}_t\|_0^2+4\mu^{-1}\|\mathbf{f}(t)-\bar{\mathbf{f}}\|_0^2\nonumber\\
&+12\mu^{-1}\bigg(\frac{\rho}{\delta^{1-\beta}}\int_{\delta t}^\infty s^{-\beta}e^{-s}ds\bigg)^2\|A\bar{\mathbf{u}}\|_{0}^2\nonumber\\
&+\frac{12\rho^2}{\mu}\frac{\Gamma(1-\beta)}{\delta^{1-\beta}}\int_{0}^{t}(t-s)^{-\beta}e^{-\delta(t-s)}\|A\mathbf{z}\|_0^2ds.\label{eq4.1}
\end{align}
Invoking (\ref{eq3.5}), (\ref{eq3.25}), (\ref{eq3.5aad}), (\ref{eq3.11ad1}) and (\ref{eq2.8}) in (\ref{eq4.1}) yields
\begin{eqnarray}
\mu \|A\mathbf{z}\|_0^2\leq \kappa (e^{-\delta t}+e^{-2\alpha t})+\frac{12\rho^2}{\mu}\frac{\Gamma(1-\beta)}{\delta^{1-\beta}}\int_{0}^{t}(t-s)^{-\beta}e^{-\delta(t-s)}\|A\mathbf{z}\|_0^2ds,\label{eq3.13adh}
\end{eqnarray}
together with $0<\alpha<\frac{1}{2}\min\{\delta,\frac{\mu_0\gamma_0}{2}\}$ and Lemma \ref{L4.2}, then (\ref{eq3.13adh}) becomes
\begin{eqnarray}
\|A\mathbf{z}\|_0^2\leq \kappa e^{-2\alpha t}\exp\bigg(\frac{12\rho^2}{\mu}\frac{\Gamma(1-\beta)}{\delta^{1-\beta}}\int_{0}^{t}(t-s)^{-\beta}e^{-\delta(t-s)}ds\bigg),\nonumber
\end{eqnarray}
which, combine with (\ref{eq2.3}) and (\ref{gambnd}), implies (\ref{eq3.15ad}).

Finally, using (\ref{eq3.7}) yields
\begin{align}
d(\mathbf{v},\eta)&=(\mathbf{z}_t,\mathbf{v})+\mu a(\mathbf{z},\mathbf{v})+a_1(\mathbf{z},\mathbf{z},\mathbf{v})+a_1(\bar{\mathbf{u}},\mathbf{z},\mathbf{v})+a_1(\mathbf{z},\bar{\mathbf{u}},\mathbf{v})\nonumber\\
&+d(\mathbf{z},q)+J(t;\mathbf{z},\mathbf{v})-J_1(t;\bar{\mathbf{u}},\mathbf{v})-(\mathbf{f}(t)-\bar{\mathbf{f}},\mathbf{v}),\ \ \forall\ \mathbf{v}\in \mathbf{V}.\label{eq3.7re}
\end{align}
Using (\ref{eq2.4}), (\ref{eq2.03}) and (\ref{eq2.04}) in (\ref{eq3.7re}) we obtain
\begin{align}
\|\eta(t)\|_1&\leq(1+\gamma_0^{-1})|\eta(t)|_1\leq C_1\sup_{\mathbf{v}\in \mathbf{X}}\frac{|d(\mathbf{v},q)|}{\|\mathbf{v}\|_0}\nonumber\\
&\leq C_1\|\mathbf{f}(t)-\bar{\mathbf{f}}\|_0+ C_1\|\mathbf{z}_t(t)\|_0+C_1|A\mathbf{z}(t)|_0\nonumber\\
&~~+C_1(\|A\mathbf{z}(t)\|_0+\|A\bar{\mathbf{u}}\|_0)|\mathbf{z}(t)|_1\nonumber\\
&~~+C_1\int_0^t(t-s)^{-\beta}e^{-\delta(t-s)}\|A\mathbf{z}\|_0ds\nonumber\\
&~~+\frac{\rho C_1}{\delta^{1-\beta}}\int_{\delta t}^\infty s^{-\beta}e^{-s}\|A\bar{\mathbf{u}}\|_0ds,\label{eq3.26t}
\end{align}
where $C_1>0$ is a positive constant. We apply (\ref{eq3.13adh}), (\ref{eq3.11ad1}) and (\ref{eq2.8}) as well as $0<\alpha<\frac{1}{2}\min\{\delta,\frac{\mu_0\gamma_0}{2}\}$ to obtain
\begin{align*}
&\int_0^t(t-s)^{-\beta}e^{-\delta(t-s)}\|A\mathbf{z}\|_0ds\nonumber\\
&\ds\qquad\leq e^{-\alpha t}\int_0^t(t-s)^{-\beta}e^{-(\delta-\alpha)(t-s)}e^{\alpha s}\|A\mathbf{z}\|_0ds\nonumber\\
&\qquad\leq e^{-\alpha t}\int_0^t(t-s)^{-\beta}e^{-(\delta-\alpha)(t-s)}ds\leq \kappa e^{-\alpha t},\\
&\frac{\rho}{\delta^{1-\beta}}\int_{\delta t}^\infty s^{-\beta}e^{-s}\|A\bar{\mathbf{u}}\|_0ds\\
&\ds\qquad\leq\kappa \frac{\rho}{\delta^{1-\beta}}\int_{\delta t}^\infty s^{-\beta}e^{-s}ds\leq\kappa e^{-\frac{\delta t}{2}}\leq\kappa e^{-\alpha t}.
\end{align*}
Combining (\ref{eq3.26t}) with the above two equations and using (\ref{eq3.25}), (\ref{eq3.5}), (\ref{eq3.13adh}) and (\ref{eq2.8}) we prove (\ref{eq3.17ad}) and thus complete the proof of the theorem.
 $\blacksquare$


%
 \section*{Conflict of interest}
 The authors declare that they have no conflict of interest.

\section*{Data availability statement}
Data sharing not applicable to this article as no datasets were generated or analysed during the current study.

\end{document}